\documentclass{ws-aupo}
\usepackage{epsfig}
\usepackage[OT2,T1]{fontenc}
 \newlanguage\fakelanguage
\newcommand\cyr{\fontencoding{OT2}\fontfamily{wncyr}\selectfont
\language\fakelanguage} \DeclareTextFontCommand{\textcyr}{\cyr}

\begin{document}

\markboth{ANDREA OSSICINI}{RIEMANN'S FUNCTIONAL EQUATION ALTERNATIVE
FORM, II}

\catchline{}{}{}{}{}

\title{AN ALTERNATIVE FORM OF THE FUNCTIONAL EQUATION\\
FOR RIEMANN'S ZETA FUNCTION, II}

\vspace{3mm}

\author{ANDREA OSSICINI}

\vspace{3mm}

\maketitle

\begin{abstract}
This paper treats about one of the most remarkable achievements {by
Riemann,that is the symmetric form of the functional equation for
$\zeta (s)$. We present here, after showing the first proof of
Riemann, a new, simple and direct proof of the symmetric form of the
functional equation for both the Eulerian Zeta function and the
alternating Zeta function, connected with odd numbers. A proof that
Euler himself could have arranged with a little step at the end of
his paper \textit{''Remarques sur un beau rapport entre les
s\'{e}ries des puissances tant direct que r\'{e}ciproches''}. This
more general functional equation gives origin to a special
function,here named \textcyr{E1}(s)\ which we prove that it can be
continued analytically to an entire function over the whole complex
plane using techniques similar to those of the second proof of
Riemann. Moreover we are able to obtain a connection between
Jacobi's imaginary transformation and an infinite series identity of
Ramanujan. Finally, after studying the analytical properties of the
function \textcyr{E1}(s),we complete and extend the proof of a
Fundamental Theorem, both on the zeros of Riemann Zeta function and
on the zeros of Dirichlet Beta function, using also the Euler-Boole
summation formula}.
\end{abstract}
\vspace{3mm}



\vspace{2mm}

\section{Introduction}

\vspace{3mm}

In [14] we introduced a special function, named $A(s)$, which is
\begin{equation}
\label{eq1} A\left( s \right)=\frac{\Gamma \left( s \right)\zeta
\left( s \right)L\left( s \right)}{\pi^{s}}\;\;\;\; with \;\;s \in
 \textbf{\textit{C}}.\end{equation}

where $\Gamma \left( s \right)$ denotes Euler's Gamma function,
$\zeta \left( s \right)$ denotes the Riemann Zeta function and
$L\left( s \right)$ denotes Dirichlet's (or Catalan's) Beta
function.

\vspace{2mm}

Let us remember that the Gamma function can be defined by the
Euler's integral of the second kind [22, p.241]:
\[
\Gamma \left( s \right) = \int_0^\infty {e^{ - t}t^{s - 1}} dt =
\int_0^1 {\left( {\log 1 \mathord{\left/ {\vphantom {1 t}} \right.
\kern-\nulldelimiterspace} t} \right)^{s - 1}dt}\quad \quad \quad
\left( {\Re \left( s \right)>0}\right)\]

\noindent and also by the following Euler's definition [22, p.237]:
\[
\Gamma \left( s \right) = \mathop {\lim }\limits_{n \to \infty }
\;\frac{1\cdot 2\cdot 3\cdot \, \cdot \, \cdot \, \left( {n -
1}\right)}{s\left( {s + 1} \right)\left( {s + 2} \right) \cdot \,
\cdot \, \cdot \left( {s + n - 1} \right)} \,\,n^{s}.
\]

\vspace{8mm}

The Riemann Zeta function is defined by ([17], pp. 96-97, see
Section 2.3) :

\vspace{3mm}

\[
\begin{array}{l}
 \zeta \left( s \right)\mbox{\, :}=\left\{ {\begin{array}{l}
 \sum\limits_{n=1}^\infty {\,\,\frac{1}{\left( n
\right)^{s}}=\frac{1}{1-2^{-s}}\sum\limits_{n=1}^\infty
{\frac{1}{\left( {2n-1} \right)^{s}}\quad \quad \quad \left( {\Re
\left( s \right)>1}
\right)} } \\
 \\
 \left( {1-2^{1-s}} \right)^{-1}\sum\limits_{n=1}^\infty {\frac{\left( {-1}
\right)^{n-1}}{\left( n \right)^{s}}\quad \quad \quad \quad \,
\left( {\Re \left( s \right)>0,s\ne 1} \right)} \\
 \end{array}} \right.\quad \quad \quad \;\quad \quad \\
 \\
 \end{array}
\]
which can be indeed analytically continued to the whole complex
\textit{ s-}plane except for a simple pole at $s=1$ with residue 1.

\vspace{3mm}

The Riemann Zeta function $\zeta \left( s \right)$ plays a central
role in the applications of complex analysis to number theory.

\vspace{3mm}

The number-theoretic properties of $\zeta \left( s \right)$ are
exhibited by the following result as \textit{Euler's product
formula}, which gives a relationship between the set of primes and
the set of positive integers:
\[
\zeta \left( s \right) = \prod\limits_p {\left( {1 - p^{ - s}}
\right)} ^{ - 1}\quad \left( {\Re \left( s \right) > 1} \right),
\]
\noindent where the product is taken over all primes.

It is an analytic version of the fundamental theorem of arithmetic,
which states that every integer can be factored into primes in an
essentially unique way.

\vspace{3mm}

Euler used this product to prove that the sum of the reciprocals of
the primes diverges.

\vspace{2mm}

The Dirichlet Beta function, also known as Dirichlet's $ L$ function
for the nontrivial character modulo 4, is defined, practically for
$\Re \left( s \right)>0$, by:
\[
L\left( s \right)=\,L\left( s,\chi_{\mbox{4}}
\right)=\sum\limits_{n=0}^\infty {\frac{\left( {-1}
\right)^{n}}{\left( {2n+1} \right)^{s}}}
\]
and it does not possess any singular point.

The $L\left( s \right)$  function is also connected to the theory of
primes which may perhaps be best summarized by

\vspace{2mm}

\[
L\left( s \right) = \prod\limits_{p{\kern 1pt} \equiv {\kern 1pt}
1{\kern 1pt} \bmod {\kern 1pt} 4} {\left( {1 - p^{ - s}} \right)} ^{
- 1} \cdot \prod\limits_{p{\kern 1pt} \equiv {\kern 1pt} 3{\kern
1pt} \bmod {\kern 1pt} 4} {\left( {1 - p^{ - s}} \right)} ^{ - 1} =
\prod\limits_{p\,\,odd{\kern 1pt} } {\left( {1 - \left( { - 1}
\right)^{\frac{p - 1}{2}}p^{ - s}} \right)} ^{ - 1},
\]

\vspace{3mm}

\noindent where the products are taken over primes and the
rearrangement of factors is permitted because of an absolute
convergence.

\vspace{6mm}

In [14] we have also proved the following identity:
\begin{equation}
\label{eq2} A\left( s \right)=A\left( {1-s} \right)\end{equation}

and we have used the functional equation of $L\left( s \right)$ to
rewrite the functional equation (\ref{eq2}) in Riemann's well known
functional equation for Zeta:
\begin{equation}
\label{eq3} \zeta \left( s \right)=2^{s}\pi^{s-1}\Gamma \left( {1-s}
\right)\sin \left( {\frac{\pi \;s}{2}} \right)\zeta \left( {1-s}
\right)\end{equation}

or equivalently to
\[\zeta \left( {1-s} \right)=2\left( {2\;\pi }
\right)^{-s}\Gamma \left( s \right)\cos \left( {\frac{\pi \;s}{2}}
\right)\zeta \left( s \right).\]

This approach is the motivation for saying that the following
symmetrical formulation:
\[\pi^{-s}\Gamma \left( s \right)\;\zeta \left( s \right)\;L\left( s
\right)=\pi ^{ - \left( {1 - s} \right)}\Gamma \left( {1-s}
\right)\;\zeta \left( {1-s} \right)\;L\left( {1-s} \right).\] is an
alternative form of the functional equation for Riemann's Zeta
Function.

\section{The origin of the symmetric form of the
functional  equation for the Eulerian Zeta and for the alternating
Zeta, connected with odd numbers.}

Riemann gives two proofs of the functional equation (\ref{eq3}) in
his paper [15], and subsequently he obtains the symmetric form by
using two basic identities of the factorial function, that are
Legendre's duplication formula [13], which was discovered in 1809
and was surely unknown to Euler:

\[\Gamma \left( s \right)\;\Gamma \left( s+\frac{1}{2}
\right)=\frac{\sqrt \pi }{2^{2s-1}}\;\Gamma \left( {2s} \right)\]

and Euler's complement formula:
\[\Gamma \left( s
\right)\;\Gamma \left( 1-s \right)=\frac{\pi }{\sin \left(\pi \;s
\right)}.\]

Riemann rewrites the functional equation (\ref{eq3}) in the form (
[6], pp. 12-15):
\[ \zeta \left( s \right)=\frac{2^{s}\pi^{s-1}}{\sqrt \pi
}2^{-s}\Gamma \left( {\frac{1-s}{2}} \right)\Gamma \left(
{1-\frac{s}{2}} \right)\frac{\pi } {\Gamma \left( {\frac{s}{2}}
\right)\Gamma \left( {1-\frac{s}{2}}  \right)}\zeta \left( {1-s}
\right) \]

\vspace{3mm}

and using the simplification $\pi^{s-1}\sqrt \pi
\;=\frac{\pi^{-{\left( {1-s} \right)} \mathord{\left/ {\vphantom
{{\left( {1-s} \right)} 2}}  \right. \kern-\nulldelimiterspace}
2}}{\pi^{-s \mathord{\left/ {\vphantom  {s 2}} \right.
\kern-\nulldelimiterspace} 2}}$, he obtains the desired formula:

\vspace{1mm}

\[ \Gamma \left( {\frac{s}{2}} \right)\;\pi^{^{-s \mathord{\left/
{\vphantom  {s 2}} \right. \kern-\nulldelimiterspace} 2}}\zeta
\left( s  \right)\;=\;\Gamma \left( {\frac{1-s}{2}}
\right)\;\pi^{^{-{\left( {1-s}  \right)} \mathord{\left/ {\vphantom
{{\left( {1-s} \right)} 2}} \right.  \kern-\nulldelimiterspace}
2}}\zeta \left( {1-s} \right).\]

\vspace{3mm}

Now this property induced Riemann to introduce, in place of $\Gamma
\left( s \right)\;$, the integral $\Gamma \left( {\frac{s}{2}}
\right)$ and at the end, for convenience, to define the $\xi $
function as:
\begin{equation}
\label{eq4} \xi \left( s \right)\;\;=\;\frac{s}{2}\left( {s-1}
\right)\pi^{-s \mathord{\left/ {\vphantom {s 2}} \right.
\kern-\nulldelimiterspace} 2}\Gamma \left( {\frac{s}{2}}
\right)\;\zeta \left( s \right).\end{equation}

In this way $\xi \left( s \right)$ is an entire function and
satisfies the simple functional equation:

\begin{equation}
\label{eq5} \xi \left( s \right)\;\;=\xi \left( {1-s}
\right)\;.\end{equation}

\vspace{2mm}

This shows that $\xi \left( s \right)\;$is symmetric around the
vertical line $\Re \left( s \right)=\frac{1}{2}$.

\vspace{3mm}

In Remark 2 of [14] we stated that Euler himself could have proved
the identity (\ref{eq2}) using three reflection formulae of the
$\zeta \left( s \right)$ ,$L\left( s \right)$ and $\Gamma \left( s
\right)$, all well-known to him.

\vspace{3mm}

Here we present the simplest and direct proof based on the
astonishing conjectures, that are Euler's main results in his work
``Remarques sur un beau rapport entre les s\'{e}ries des puissances
tant direct que r\'{e}ciproches'' [8].

\vspace{2mm}

Euler writes the following functional equations:

\[\frac{1-2^{n-1}+3^{n-1}-4^{n-1}+5^{n-1}-6^{n-1}+\;\cdot
\cdot \cdot }{1-2^{-n}+3^{-n}-4^{-n}+5^{-n}-6^{-n}+\cdot \cdot \cdot
}=-\frac{1\cdot 2\cdot 3\cdot \cdot \cdot \left( {n-1}
\right)\,\left( {2^{n}-1} \right)}{\left( {2^{n-1}-1}
\right)\,\,\pi^{n}}\cos \left( {\frac{n\pi }{2}} \right)\]

and

\[\frac{1-3^{n-1}+5^{n-1}-7^{n-1}+\;\cdot \cdot \cdot
}{1-3^{-n}+5^{-n}-7^{-n}+\cdot \cdot \cdot }=\frac{1\cdot 2\cdot
3\cdot \cdot \cdot \left( {n-1} \right)\,\left( {2^{n}}
\right)}{\pi^{n}}\sin \left( {\frac{n\pi }{2}} \right)_{}\]

and concludes his work by proving that those conjectures are valid
for positive and negative integral values as well as for fractional
values of $n$.

\vspace{2mm}

In modern notation we have, with $s\in
 \textbf{\textit{C}}$:

\begin{equation}
\label{eq6} \frac{\eta \left( {1-s} \right)}{\eta \left( s
\right)}=-\frac{\left( {2^{s}-1} \right)}{\pi^{s}\left( {2^{s-1}-1}
\right)}\Gamma \left( s \right)\cos \left( {\frac{s\pi }{2}}
\right)\end{equation}

and

\begin{equation}
\label{eq7} \frac{L\left( {1-s} \right)}{L\left( s
\right)}=\frac{2^{s}}{\pi^{s}}\Gamma \left( s \right)\sin \left(
{\frac{s\pi }{2}} \right).\end{equation}

(\ref{eq6}) represents the functional equation of Dirichlet's Eta
function, which is defined for $\Re \left( s \right)\succ 0$ through
the following alternating series:

\[\eta \left( s
\right)=\sum\limits_{n=1}^\infty {\frac{\left( {-1}
\right)^{n-1}}{n^{s}}} .\]

This function $\eta \left( s \right)$ is one simple step removed
form $\zeta \left( s \right)$ as shown by the relation:

\[\eta \left( s \right)=\quad \left( {1-2^{1-s}} \right)\,\zeta
\left( s \right).\]

\vspace{1mm}

Thus (\ref{eq6}) is easily manipulated into relation (\ref{eq3}).

\vspace{1mm}

The (\ref{eq7}) is the functional equation of Dirichlet's $L$
function.

\vspace{1mm}

That being stated, multiplying (\ref{eq6}) by (\ref{eq7}) we obtain:

\[\frac{\eta
\left( 1-s \right)}{\eta \left( s \right)}\cdot \frac{L\left( 1-s
\right)}{L\left( s \right)}=\frac{\left( {1-2^{s}} \right)\cdot
2^{s-1}\left[ {\Gamma \left( s \right)} \right]^{2}}{\pi^{s}\left(
{2^{s-1}-1} \right)\,\pi^{s-1}}\cdot \frac{2\sin \left( {\frac{s\pi
}{2}} \right)\cos \left( {\frac{s\pi }{2}} \right)}{\pi }.\]

\vspace{1mm}

Considering the duplication formula of $\sin \left( s \pi \right)$
and Euler's complement formula we have:

\[\frac{\eta \left( {1-s} \right)}{\eta \left( s \right)}\cdot
\frac{L\left( 1-s \right)}{L\left( s \right)}=\frac{\left( {1-2^{s}}
\right)\,\,\pi ^{1-s}}{\left( {1-2^{1-s}} \right)\,\,\pi^{s}}\cdot
\frac{\Gamma \left( s \right)}{\Gamma \left( {1-s} \right)}.\]

\vspace{1mm}

Shortly and ordering we obtain the following remarkable identity:

\begin{equation}
\label{eq8} \frac{\left( {1-2^{1-s}} \right)}{\pi^{1-s}}\cdot \Gamma
\left( {1-s} \right)\,\eta \left( {1-s} \right)\,L\left( {1-s}
\right)=\frac{\left( {1-2^{s}} \right)}{\pi^{s}}\cdot \Gamma \left(
s \right)\,\eta \left( s \right)\,L\left( s \right).\end{equation}

\vspace{1mm}

This is unaltered by replacing $\left( {1-s} \right)$ by $s$.

\vspace{1mm}

\section{The special function $\textcyr{E1}\left( s \right)$ and its
integral representation}

\vspace{1mm}

At this stage, let us introduce the following special
function\footnote{ The letter \textcyr{E1}, called E reversed, is a
letter of the Cyrillic alphabet and is the third last letter of the
Russian alphabet.}:

\begin{equation}
\label{eq9} \textcyr{E1}\left( s \right)=\frac{\left( {1-2^{s}}
\right)\,\Gamma \left( s \right)\,\eta \left( s \right)\,L\left( s
\right)}{\pi^{s}}=\frac{\left( {1-2^{s}} \right)\left( {1-2^{1-s}}
\right)\,\Gamma \left( s \right)\zeta \left( s \right)L\left( s
\right)}{\pi^{s}}.\end{equation}

\vspace{1mm}

It is evident that from (\ref{eq1}) one has:\[\textcyr{E1}\left( s
\right)=\left( {1-2^{s}} \right)\left( {1-2^{1-s}} \right)\,A\left(
s \right).\]

\vspace{1mm}

This choice is based upon the fact that $\textcyr{E1}\left( s
\right)$ is an entire function of $s$, hence it has no poles and
satisfies the simple functional equation:
\begin{equation}
\label{eq10} \textcyr{E1}\left( s \right)=\textcyr{E1}\left( {1-s}
\right).\end{equation}

The poles at $s=0,1$ , respectively determined by the Gamma function
$\Gamma \left( s \right)$ and by the Zeta function $\zeta \left( s
\right)$ are  cancelled by the term $\left( {1-2^{s}} \right)\cdot
\left( {1-2^{1-s}} \right)$.

\vspace{1mm}

Now by using the identities ([5], chap. X, p. 355,10.15):

\begin{equation}
\label{eq11} \Gamma \left( s \right)\,a^{-s}=\int_0^\infty
{x^{s-1}e^{-ax}dx\equiv M_s} \left\{ {\,e^{-ax}\,}
\right\}\end{equation}

where $M_{s} $, denotes the Mellin transform and
\begin{equation}
\label{eq12} \sum\limits_m e^{ - m^2x} = \frac{1}{2}\left[ {\theta
_3 \left( {0\left| {{ix} \mathord{\left/ {\vphantom {{ix} \pi }}
\right. \kern-\nulldelimiterspace} \pi } \right.} \right) - 1}
\right]
\end{equation}
where $\theta_{3} \left( {z\left| \tau \right.} \right)$ is one of
the four theta functions, introduced by of Jacobi ([22], chap. XXI)
and the summation variable $m$ is to run over all positive integers,
we derive the following integral representation of $\textcyr{E1}
\left( s \right)$ function:
\begin{equation}
\label{eq13} \textcyr{E1}\left( s \right)=\frac{\left( {1-2^{s}}
\right)}{2}\frac{\left( {1-2^{1-s}}
\right)}{2}\,\;\int\limits_0^\infty \left[ {\theta_{3}^{2} \left(
{0\left| {{ix} \mathord{\left/ {\vphantom {{ix} \pi }} \right.
\kern-\nulldelimiterspace} \pi } \right.} \right) - 1} \right] \cdot
\left( {\frac{x}{\pi }} \right)^{s}\frac{dx}{x}\end{equation} Indeed
combining the following two Mellin transforms:
\[\Gamma \left( s \right)\zeta \left( {2s} \right)=M_{s} \left\{
{\frac{1}{2}\left[ {\,\theta_{3} \left( {0\left| {{ix}
\mathord{\left/ {\vphantom {{ix} \pi }} \right.
\kern-\nulldelimiterspace} \pi } \right.} \right)-1} \right]}
\right\}\;\;\;\;\left( {\Re \left( s \right)>\frac{1}{2}} \right)\]
and
\[\Gamma \left( s \right)\left[ {L\left( s \right)\zeta \left( s
\right)-\zeta \left( {2s} \right)} \right]=M_{s} \left\{
{\frac{1}{4}\left[ {\,\theta_{3} \left( {0\left| {{ix}
\mathord{\left/ {\vphantom {{ix} \pi }} \right.
\kern-\nulldelimiterspace} \pi } \right.} \right)-1} \right]\;^{2}}
\right\}\quad\left( {\Re \left( s \right)>1} \right)_{,}\] the
former is immediately obtained from Eqs. (\ref{eq11}) and
(\ref{eq12}) and the latter is obtained integrating term by term the
following remarkable identity, obtained from an identity by Jacobi
[11] and the result\footnote{ $K$ denotes the complete elliptic
integral of the first kind of modulus $k$\par } $\theta_{3}^{2}
\left( {0\left| \tau \right.} \right)={2K} \mathord{\left/
{\vphantom {{2K} \pi }} \right. \kern-\nulldelimiterspace} \pi $
([22], p. 479):
\[\frac{1}{4}\left[ {\theta_{3}^{2}
\left( {0\left| {{ix} \mathord{\left/ {\vphantom {{ix} \pi }}
\right. \kern-\nulldelimiterspace} \pi } \right.} \right)-1}
\right]\;=\sum\limits_\ell {\left( {-1} \right)^{{\left( {\ell -1}
\right)} \mathord{\left/ {\vphantom {{\left( {\ell -1} \right)} 2}}
\right. \kern-\nulldelimiterspace} 2}} \left[ {e^{\ell x}-1}
\right]^{-1}\](here the sum is to expanded as a geometric series in
$e^{-\ell x}$:\[e^{-\ell x}+e^{-2\ell x}+e^{-3\ell x}+e^{-4\ell
x}+e^{-5\ell x}+\;\cdot \cdot \cdot =\left[ {e^{\ell x}-1}
\right]^{-1}\]and the summation variable $\ell $ runs over all
positive odd integers), thus we are in the position to determine the
integral representation (\ref{eq13}) for the $\textcyr{E1}\left( s
\right)$ function, by the linearity property of Mellin
transformation, from the following identity:
\begin{equation}
\label{eq14} \textcyr{E1}\left( s \right)=\frac{\left( {1-2^{s}}
\right)\,\left( {1-2^{1-s}} \right)\,\Gamma \left( s
\right)\,\,\zeta \left( s \right)\,L\left( s
\right)}{\pi^{s}}\end{equation}
\[=\frac{\left( {1-2^{s}} \right)\left( {1-2^{1-s}}
\right)}{\pi^{s}}\,\cdot \left[ {M_{s} \left\{ {\frac{1}{4}\left[
{\,\theta_{3} \left( {0\left| {{ix} \mathord{\left/ {\vphantom {{ix}
\pi }} \right. \kern-\nulldelimiterspace} \pi } \right.} \right)-1}
\right]\;^{2}} \right\}-M_{s} \left\{ {\frac{1}{2}\left[
{\,\theta_{3} \left( {0\left| {{ix} \mathord{\left/ {\vphantom {{ix}
\pi }} \right. \kern-\nulldelimiterspace} \pi } \right.} \right)-1}
\right]} \right\}} \right]\]
\[=\frac{\left( {1-2^{s}} \right)\cdot \left( {1-2^{1-s}} \right)}{\pi^{s}}
\;M_{s} \left\{ {\frac{1}{4}\left[ {\theta_{3}^{2} \left( {0\left|
{{ix} \mathord{\left/ {\vphantom {{ix} \pi }} \right.
\kern-\nulldelimiterspace} \pi } \right.} \right)-1}
\right]}\right\} \; \quad\left( {\Re \left( s \right)>1} \right).\]

\vspace{1mm}

Now we start from (\ref{eq13}) to give an independent proof of
(\ref{eq10}) that does not use (\ref{eq6}) and (\ref{eq7}), adopting
techniques similar to Riemann's ones we use the following
fundamental transformation formula for $\theta _3 \left( {z\left|
\tau \right.} \right)$:
\begin{equation}
\label{eq15} \theta_{3} \left( {z\left| \tau \right.} \right)=\left(
{-i\tau }  \right)^{-1 \mathord{\left/ {\vphantom {1 2}} \right.
\kern-\nulldelimiterspace} 2}\exp \,\left( {{z^{2}} \mathord{\left/
{\vphantom {{z^{2}} {\pi \;i\,\tau }}} \right.
\kern-\nulldelimiterspace}  {\pi \;i\,\tau }} \right)\cdot
\theta_{3} \left( {\frac{z}{\tau } \left| {-\frac{1}{\tau }}
\right.} \right) \end{equation}

where $\left({- i\tau } \right)^{-1 \mathord{\left/{\vphantom {1
2}}\right. \kern-\nulldelimiterspace} 2}$ is to be interpreted by
the convention  $\left|{\arg \left({-i\tau
}\right)}\right|<\frac{1}{2}\pi$ ([22],p. 475).

In particular we obtain that:

\begin{equation}
\label{eq16} \theta_{3}^{2} \left( {0\left| {\frac{i\,x}{\pi }}
\right.} \right)=\frac{\pi }{x}\theta_{3}^{2} \left( {0\left|
{\frac{i\,\pi \,}{x}} \right.} \right). \end{equation}

We then rewrite the integral that appears in (\ref{eq13}) as:
\[
\int\limits_0^\pi {\left[ {\theta_{3}^{2} \left( {0\left| {{ix}
\mathord{\left/ {\vphantom {{ix} \pi }} \right.
\kern-\nulldelimiterspace} \pi } \right.} \right)-1} \right]} \cdot
\left( {\frac{x}{\pi }}
\right)^{s}\frac{dx}{x}+\int\limits_\pi^\infty {\left[
{\theta_{3}^{2} \left( {0\left| {{ix} \mathord{\left/ {\vphantom
{{ix} \pi }} \right. \kern-\nulldelimiterspace} \pi } \right.}
\right)-1} \right]} \cdot \left( {\frac{x}{\pi }}
\right)^{s}\frac{dx}{x}
\]
\[
=\int\limits_0^\pi {\left[ {\theta_{3}^{2} \left( {0\left| {{ix}
\mathord{\left/ {\vphantom {{ix} \pi }} \right.
\kern-\nulldelimiterspace} \pi } \right.} \right)} \right]} \cdot
\left( {\frac{x}{\pi }}
\right)^{s}\frac{dx}{x}-\frac{1}{s}+\int\limits_\pi^\infty {\left\{
{\left[ {\theta_{3}^{2} \left( {0\left| {{ix} \mathord{\left/
{\vphantom {{ix} \pi }} \right. \kern-\nulldelimiterspace} \pi }
\right.} \right)-1} \right]} \right\}} \cdot \left( {\frac{x}{\pi }}
\right)^{s}\frac{dx}{x}
\]

and use the change of variable $\frac{i\,x}{\pi }\to \frac{i\,\pi
}{x}$ and the (\ref{eq16}) to find:

\vspace{1mm}

\[
\int\limits_0^\pi {\;\left[ {\theta_{3}^{2} \left( {0\left| {{ix}
\mathord{\left/ {\vphantom {{ix} \pi }} \right.
\kern-\nulldelimiterspace} \pi } \right.} \right)} \right]} \cdot
\left( {\frac{x}{\pi }} \right)^{s}\frac{dx}{x}=\int\limits_{\;\pi
}^\infty {\;\left[ {\theta _{3}^{2} \left( {0\left| {{i\,\pi }
\mathord{\left/ {\vphantom {{i\,\pi } x}} \right.
\kern-\nulldelimiterspace} x} \right.} \right)} \right]} \cdot
\left( {\frac{\pi }{x}}
\right)^{s}\frac{dx}{x}=\int\limits_\pi^\infty {\;\left[
{\theta_{3}^{2} \left( {0\left| {{i\,x} \mathord{\left/ {\vphantom
{{i\,x} \pi }} \right. \kern-\nulldelimiterspace} \pi } \right.}
\right)} \right]} \cdot \left( {\frac{x}{\pi }}
\right)^{1-s}\frac{dx}{x}
\]
\[
=-\frac{1}{1-s}+\int\limits_\pi^\infty {\;\left[ {\theta_{3}^{2}
\left( {0\left| {{i\,x} \mathord{\left/ {\vphantom {{i\,x} \pi }}
\right. \kern-\nulldelimiterspace} \pi } \right.} \right)-1}
\right]} \cdot \left( {\frac{x}{\pi }} \right)^{1-s}\frac{dx}{x}.
\]

Therefore:

\begin{multline}
\label{eq17} \textcyr{E1}\left( s \right)=\frac{\left( {1-2^{s}}
\right)}{2}\frac{\left( {1-2^{1-s}} \right)}{2}\cdot\\
\cdot\left\{{\frac{1}{s\left( {s-1}\right)}+\int\limits_\pi^\infty
{\;\left[ {\theta_{3}^{2} \left( {0\left| {{ix} \mathord{\left/
{\vphantom {{ix} \pi }} \right. \kern-\nulldelimiterspace} \pi }
\right.} \right)-1} \right]}\cdot\left[ {\left( {\frac{x}{\pi }}
\right)^{s}+\left( {\frac{x}{\pi }} \right)^{1-s}} \right]d\log x}
\right\}\end{multline}

which is manifestly symmetrical under $s\to 1-s$, and analytic since
$\theta _{3} \left( {0\left| {\frac{i\;x}{\pi }} \right.} \right)$
decreases exponentially as $x\to \infty $.

This concludes the proof of the functional equation and the analytic
continuation of $\textcyr{E1}\left( s \right)$\textbf{, }assuming
the identity (\ref{eq15}), due to Jacobi.

\section{Jacobi's imaginary transformation and an infinite series identity of Ramanujan}

The fundamental transformation formula of Jacobi for $\theta_{3}
\left( {z\left| \tau \right.} \right)$:
\[ \theta_{3} \left(
{z\left| \tau \right.} \right)=\left( {-i\tau }  \right)^{-1
\mathord{\left/ {\vphantom {1 2}} \right. \kern-\nulldelimiterspace}
2}\exp \,\left( {{z^{2}} \mathord{\left/ {\vphantom {{z^{2}} {\pi
\;i\,\tau }}} \right. \kern-\nulldelimiterspace}  {\pi \;i\,\tau }}
\right)\cdot \theta_{3} \left( {\frac{z}{\tau } \left|
{-\frac{1}{\tau }} \right.} \right) \]

\noindent where the squares root is to be interpreted as the
principal value; that  is,$\qquad\qquad$  if $w=re^{i\theta }$ where
$0 \le
 \theta \le 2\pi $, then $w^{1 \mathord{\left/ {\vphantom {1 2}}  \right.
\kern-\nulldelimiterspace} 2}=r^{1 \mathord{\left/ {\vphantom {1 2}}
\right. \kern-\nulldelimiterspace} 2}e^{{i\theta } \mathord{\left/
{\vphantom {1 2}} \right. \kern-\nulldelimiterspace} 2}$ and the
infinite series identity of Ramanujan ([3], Entry 11, p. 258):
\[
\alpha \left\{ {\frac{1}{4}\sec \left( {\alpha n} \right) +
\sum\limits_{k = 1}^\infty {\;\chi \left( k \right)\frac{\cos \left(
{\alpha \,nk} \right)}{e^{\alpha ^2k} - 1}} } \right\} = \beta
\left\{ {\frac{1}{4} +{\frac{1}{2}} \sum\limits_{k = 1}^\infty
{\;\frac{\cosh \left( {2\beta \,nk} \right)}{\cosh \left( {\beta
^2k} \right)}} } \right\}
\]

\vspace{1mm}

\noindent\textit{with }$ \alpha ,\beta \succ 0$\textit{,}$\; \alpha
\beta = \pi $\textit{,} $n$ $\in \Re$\textit{,} $\left| n \right|
\prec \beta \mathord{\left/ {\vphantom {\beta 2}} \right.
\kern-\nulldelimiterspace} 2\;$\textit{and with} \[\;\chi \left( k
\right) = \left\{ {\begin{array}{l}
 \,0\,\quad \;for\;k\;even \\
 \;1\,\quad for\;k\; \equiv \;1\;\bmod \;4 \\
 -1\,\,\,\,for\;k\; \equiv \;3\;\bmod \;4 \\
 \end{array}} \right.
\]

\vspace{1mm}

\noindent can be derived from the following Poisson summation
formula ( see [2] , pp. 7-11 and [14], Appendix):
\[
\frac{1}{2}f\left( 0 \right)+\sum\limits_{k=1}^\infty {f\left( k
\right)} =\,\;\int\limits_0^\infty {f\left( x \right)dx+2}
\sum\limits_{k=1}^\infty {\;\int\limits_0^\infty {f\left( x \right)}
} \cos \left( {2k\pi \;x} \right)dx.
\]

\vspace{1mm}

From Jacobi's Lambert series formula for $\theta_{3}^{2} \left(
{0\left| \tau \right.} \right) $:

\[
\theta_{3}^{2} \left( {0\left| \tau \right.}
\right)-1=4\sum\limits_\ell {\left( {-1} \right)^{{\left( {\ell -1}
\right)} \mathord{\left/ {\vphantom {{\left( {\ell -1} \right)} 2}}
\right. \kern-\nulldelimiterspace} 2}} \,q^{\ell }\left( {1-q^{\ell
}} \right)^{-1},
\]
where $\ell $ is to run over all positive odd integers, we have
again with $q=\exp \left( {i\pi \,\tau } \right)_{\, }$ and $\tau
={ix} \mathord{\left/ {\vphantom {{ix} \pi }} \right.
\kern-\nulldelimiterspace} \pi $:
\[
\frac{1}{4}\left[ {\theta_{3}^{2} \left( {0\left| {{ix}
\mathord{\left/ {\vphantom {{ix} \pi }} \right.
\kern-\nulldelimiterspace} \pi } \right.} \right)-1}
\right]\;=\sum\limits_\ell {\left( {-1} \right)^{{\left( {\ell -1}
\right)} \mathord{\left/ {\vphantom {{\left( {\ell -1} \right)} 2}}
\right. \kern-\nulldelimiterspace} 2}} \left[ {e^{\ell x}-1}
\right]^{-1}.\]

Now
\[
\sum\limits_\ell {\left( {-1} \right)^{{\left( {\ell -1} \right)}
\mathord{\left/ {\vphantom {{\left( {\ell -1} \right)} 2}} \right.
\kern-\nulldelimiterspace} 2}} \left[ {e^{\ell x}-1}
\right]^{-1}=\sum\limits_{m=1}^\infty {\chi \left( m \right)}
\frac{1}{e^{mx}-1}
\]
where still
\[
\chi \left( m \right)=\left\{ {\begin{array}{l}
 0\;\quad \;for\,m\;\;even \\
 1\quad \,\,\;for\,m\equiv 1\,\,\bmod \,\,4 \\
 -\;1\,\,{\kern 1pt}\;{\kern 1pt}{\kern 1pt}for\,m\equiv 3\,\,\bmod \,\,4 \\
 \end{array}} \right.\]
and therefore
\begin{equation}
\label{eq18} \frac{1}{4}\left[ {\theta_{3}^{2} \left( {0\left| {{ix}
\mathord{\left/ {\vphantom {{ix} \pi }} \right.
\kern-\nulldelimiterspace} \pi } \right.} \right)-1}
\right]\;=\sum\limits_{m=1}^\infty {\chi \left( m \right)}
\frac{1}{e^{mx}-1}\end{equation}

For $n = 0$ the infinite series identity of Ramanujan reads

\[\alpha \left\{ {\frac{1}{4} +
\sum\limits_{k = 1}^\infty {\;\chi \left( k
\right)\frac{1}{e^{\alpha ^2k} - 1}} } \right\} = \beta \left\{
{\frac{1}{4} +{\frac{1}{2}} \sum\limits_{k = 1}^\infty
{\;\frac{1}{\cosh \left( {\beta ^2k} \right)}} } \right\}.
\]

Replacing $\cosh \left( x \right)$ by the exponential functions,
expanding the geometric series and rearranging the sums we have

\[ {\frac{1}{2}}\sum\limits_{k = 1}^\infty
{\;\frac{1}{\cosh \left( {\beta ^2k} \right)}} = \sum\limits_{m =
1}^\infty {\;\chi \left( m \right)\frac{1}{e^{\beta ^2m} - 1}}.
\]

Now we substitute $\alpha = \sqrt x $, $\beta = \pi \mathord{\left/
{\vphantom {\pi {\sqrt x }}} \right. \kern-\nulldelimiterspace}
{\sqrt x }_{ }$ and we obtain:
\[
\left\{ {\frac{1}{4} + \sum\limits_{k = 1}^\infty {\chi \left( k
\right)\frac{1}{e^{xk} - 1}} } \right\} = \frac{\pi }{x}\left\{
{\frac{1}{4} + \sum\limits_{m = 1}^\infty {\chi \left( m
\right)\frac{1}{e^{\left(\pi^2 \mathord{\left/ {\vphantom {\pi x}}
\right. \kern-\nulldelimiterspace} x\right)m} - 1}} } \right\}.\]

\vspace{1mm}

Finally, with the relation (\ref{eq18}) we establish the following
transformation of $\theta_{3}^{2} \left( {0\left| {\frac{i\,x}{\pi
}} \right.} \right)$:
\[\theta_{3}^{2} \left( {0\left| {\frac{i\,x}{\pi }}
\right.} \right)=\frac{\pi }{x}\theta_{3}^{2} \left( {0\left|
{\frac{i\,\pi \,}{x}} \right.} \right).\]

This last transformation is also an immediate consequence of the
fundamental transformation formula of Jacobi for $\theta_{3} \left(
{z\left| \tau \right.} \right)$.

In this way we have obtained an amazing connection between the
Jacobi imaginary transformation and the infinite series identity of
Ramanujan.

\vspace{1mm}

\section{The properties of the function $\textcyr{E1} \left( s
\right)$}

\vspace{2mm}

In this section we remark the following fundamental properties of
the special function $\textcyr{E1}\left( s \right)$ with $s=\sigma
+it$\textbf{:}

\vspace{3mm}

(a) $\textcyr{E1}\left( s \right)=\textcyr{E1}\left( {1-s} \right)$

\vspace{3mm}

(b) $\textcyr{E1}\left( s \right)$ is an entire function and  $
\textcyr{E1}\left( s \right) = \overline {\textcyr{E1}\left(
{\overline s } \right)}
 $

\vspace{3mm}

(c) $\textcyr{E1} \left( {\frac{1}{2}+it} \right)\in \Re $

\vspace{3mm}

(d) $\textcyr{E1} \left( 0 \right)=\textcyr{E1} \left( 1
\right)=-\frac{\log 2}{4}$

\vspace{3mm}

(e) if $\textcyr{E1}\left( s \right)=0$, then $0\le \sigma \le 1$

\vspace{3mm}

(f) $\textcyr{E1}\left( s \right)<0$ for all $s\in \Re $\textbf{.}

\vspace{5mm}

Outline of proof:

\vspace{1mm}

Using the topics developed at the end of Sections 2 and 3, the
functional equation (a) follows.

\vspace{1mm}

Regarding (b), the second expression in the definition (3.1) shows
at once that $\textcyr{E1}\left( s \right)$ is holomorphic for
$\sigma \ge 0$\textbf{, }since the simple pole of $\Gamma \left( s
\right)$ at $s=0$ and the simple pole of $\zeta \left( s \right)$ at
$s=1$ are removed by the factors $\left( {1-2^{s}} \right)$ and
$\left( {1-2^{1-s}} \right)$, and there are no poles for $\sigma \ge
0$\textbf{,} but the (a) implies $\textcyr{E1}\left( s \right)$
holomorphic on all C.

\vspace{1mm}

The second part of (b) follows from the fact that
$\textcyr{E1}\left( s \right)$ is real on the real line, thus
$\textcyr{E1}\left( s \right)-\overline {\textcyr{E1}\left(
{\overline s } \right)}_{\mathbf{\, }}$is an analytic function
vanishing on the real line, hence zero since the zeros of an
analytic function which is not identically zero can have no
accumulation point.

\vspace{1mm}

We note that $s=\frac{1}{2}+it$ where $t$ is real, then $\overline s
$ and $1-s$ coincide, so this implies (c).

\vspace{1mm}

The known values $L\left( 1 \right)=\frac{\pi }{4}$ , $\eta \left( 1
\right)=\log 2\;$ and $\;\lim\limits_{s\to 1} \frac{\left( {1-2^{s}}
\right)}{\pi^{s}}\;\Gamma \left( s \right)=-\frac{1}{\pi }$ imply
(d) for $\textcyr{E1}\left( 1 \right)$ and the functional equation
(a) then gives the result for $\textcyr{E1}\left( 0 \right)$.

\vspace{1mm}

Since the Gamma function has no zeros and since the Dirichlet Beta
function and the Riemann Zeta function have respectively an Euler
product ( see \S 1. Introduction or [9], p. 53 and p. 40):

\[
L\left( s \right)=\prod\limits_{\begin{array}{l}
 \quad \quad p\; \\ prime\;odd \\ \end{array}} {\left( {1-\left( {-1}
\right)^{\frac{p-1}{2}}\cdot p^{-s}} \right)}^{-1}\; ;\qquad \qquad
\zeta \left( s \right)=\prod\limits_{\begin{array}{l}
 \quad \;p \\ prime\; \\\end{array}} {\left( {1-p^{-s}} \right)}^{-1}
\]

\vspace{1mm}

\noindent which shows that they are non-vanishing in the right half
plane $\Re \left( s \right)>1$, the function $\textcyr{E1}\left( s
\right)$ has no zeros in $\Re \left( s \right)>1_{\, }$ and by
functional equation (a), it also has non zeros in $\Re \left( s
\right)<0$.

\vspace{1mm}

Thus all the zeros have their real parts between 0 and 1 (including
the extremes) and this proves (e).

\vspace{1mm}

Finally, to prove (f) first we note from the following integral
representations $ $ ([7], p. 1, p. 32 and p. 35) :
\[
\Gamma \left( s \right)=\int\limits_0^\infty
{\;\frac{x^{s-1}}{e^{x}}} dx\; ; \; L\left( s
\right)=\frac{1}{\Gamma \left( s \right)}\int\limits_0^\infty
{\;\frac{x^{s-1}}{e^{x}+e^{-x}}} dx\; ; \; \eta \left( s
\right)=\frac{1}{\Gamma \left( s \right)}\int\limits_0^\infty
{\;\frac{x^{s-1}}{e^{x}+1}} dx \; \left( {\Re \left( s \right)>0}
\right)
\]
that $\Gamma \left( s \right)$, $L\left( s \right)$, $\eta \left( s
\right)$ are positives for all $s\in \Re ,\;s>0$.

\vspace{3mm}

Then combining this with the negative factor
$\frac{1-2^{s}}{\pi^{s}}$ for $s>0$ the definition (\ref{eq9})
proves (f) for $s>0,\;s\ne 0$ and combining this with (d) then it
proves (f) for $s\ge 0$, whence the functional equation (a) shows
that (f) holds for all $s\in \Re $.

\vspace{8mm}

\section{The zeros of the entire function $\textcyr{E1}\left( s \right)$ and an estimate for the
number of these in the critical strip $0\le \sigma \le 1$}

\vspace{3mm}

We summarized and extended the results of the previous section in
the following theorem:

 \vspace{3mm}

\textbf{Fundamental Theorem}: (i) \textit{The zeros of
}$\textcyr{E1}\left( s \right)$ \textit{ (if any exits) are all
situated in the strip }$0\le \sigma \le 1$\textit{ and lie
symmetrically about the lines }$t=0$\textit{ and }$\sigma
=\frac{1}{2}$.

\vspace{2mm}

(ii) \textit{The zeros of} $\textcyr{E1}\left( s \right)$
\textit{are identical to the imaginary zeros of the factor} $\left(
{1-2^{s}} \right)\cdot \left( {1-2^{1-s}} \right)$ \textit{and to
the non-trivial zeros of the functions }$L\left( s
\right)$\textit{and }$\zeta \left( s \right)$; $\textcyr{E1}\left( s
\right)$ \textit{has no zeros on the real axis.}

\vspace{2mm}

(iii) T\textit{he number} $N\left( T \right)$ \textit{of zeros of}
$\textcyr{E1}\left( s \right)$ \textit{in the rectangle with}
\textbf{ }$0\le \sigma \le 1$, $0\le t\le T$,\textit{ when } $T\to
\infty$ \textit{satisfies} :
\[
N\left( T \right)=\frac{T}{\pi }\log \frac{2T}{\pi \, e}+O\left(
{\log T} \right)
\]

\textit{where the notation }$f\left( T \right)=O\left( {g\left( T
\right)} \right)$\textit{ means } $\frac{f\left( T \right)}{g\left(
T \right)}$ \textit{is bounded by a constant independent of } $T$.

\vspace{2mm}

\textbf{Proof.} To prove (i) the properties (a) and (e) are
sufficient.

\vspace{2mm}

These properties together with (b) and (c) show that we may detect
zeros of $\textcyr{E1}\left( s \right)$ on the line $\sigma =
\frac{1}{2}$ by detecting sign changes, for example, in
$\textcyr{E1} \left( {\frac{1}{2}+it} \right)$, so it is not
necessary to compute exactly the location of a zero in order to
confirm that it is on this line.

\vspace{1mm}

Thus we compute
\[
\textcyr{E1}\left( {\frac{1}{2}+5i} \right)=-\, 2.519281933...\cdot
10^{-\, 3} ;\quad \textcyr{E1}\left( {\frac{1}{2}+7i} \right)=+\,
8.959203701...\cdot 10^{-\, 5}
\]
we know that there is a zero of $\textcyr{E1}\left( {\frac{1}{2}+it}
\right)$ with $t$ between $5$ and $7$\textbf{.}

\vspace{1mm}

Indeed for $t=\mbox{6.0209489...}\,$\textbf{, }we have $
\textcyr{E1}\left( {\frac{1}{2}+6.0209489...i} \right)=0$ and this
is the smallest zero of $\textcyr{E1}\left( s \right)$: a much
smaller value than the one corresponding to $\zeta \left( s
\right)$, that is $\zeta \left({\frac{1}{2}+14.13472514...i }
\right)=0.$

\vspace{1mm}

To prove (ii) we have:
\[
\textcyr{E1}\left( s \right)=h\left( s \right)\frac{\Gamma \left( s
\right)\zeta \left( s \right)L\left( s \right)}{\pi^{s}}
\]

\vspace{1mm}

where the imaginary zeros of the factor $h\left( s \right)=\left(
{1-2^{s}} \right)\cdot \left( {1-2^{1-s}} \right)$ lie on the
vertical lines $\Re \left( s \right)=0$ and $\Re \left( s
\right)=1$.

\vspace{1mm}

We recall the following identity of the general exponential function
$w=a^{z}$ ( $a\ne 0$ is any complex number): $a^{z}=e^{z\log a}$;
now, the function $e^{z}$ assumes all values except zero, i.e. the
equation $e^{z}=A$ is solvable for any nonzero complex number $A$.

\vspace{1mm}

If $\alpha =\arg A$, all solutions of the equation $e^{z}=A$ are
given by the formula:
\[
z=\log \left| A \right|+i\left( {\alpha +2k\pi } \right),\quad \quad
\quad k=0,\;\pm 1,\;\pm 2,\;...
\]

In particular, if $e^{z}=1$, we have $z=2k\pi \;i,\quad k=0,\;\pm
1,\;\pm 2,\;...$.

\vspace{1mm}

Consequently, the imaginary roots of $\;h\left( s \right)\;$ are
$s=\pm \frac{2\pi \,i\,k}{\log 2}\;$ and $\; s=1\pm \frac{2\pi
\,i\,k}{\log 2}\quad $ with $k\in N,\;k>0$.

In addition from each of functional equations (2.3) and (2.4),
exploiting the zeros of the trigonometric function \textit{cosine}
and \textit{sine}, it is immediate to verify that:
\[\zeta \left( s \right)=0 \quad for \quad s=-2,\;-4,\;-6,\;-8,\cdots\]
and
\[L\left( s \right)=0 \quad for \quad s=-1,\;-3,\;-5,\;-7,\cdots .\]
These are the trivial zeros of the two Euler's Zeta functions $\zeta
\left( s \right)$ and $L\left( s \right)$, that are cancelled by the
singularities of the $\Gamma \left( s \right)$ function in the
negative horizontal axis $x.$

We remember that the two last singularities at $s=0,\;1$,
respectively determined by the $\Gamma \left( s \right)$ function
and by $\zeta \left( s \right)$ function, are cancelled by real
roots of factor $h\left( s \right)$.

We've still got the non-trivial zeros of the functions $\zeta \left(
s \right)$ and $L\left( s \right)$, see Section 5 and at the end
let's see also the property (f).

For the proof of (iii) we consider the fact that $\textcyr{E1}\left(
s \right)$ is an entire function of $s$, hence it has no poles and
the result (ii).

These properties can be then used to estimate $N\left( T \right)$ by
calling upon the Argument Principle ( [10], pp. 68-70).

\vspace{3mm}

The Argument Principle is the following theorem of Cauchy:

\vspace{2mm}

\textbf{Theorem 6.1.} \textit{Suppose the function }$F\left( s
\right)$\textit{ is analytic, apart from a finite number of poles,
in the closure of a domain D bounded by a simple closed positively
oriented Jordan curve C. Suppose further that }$F\left( s
\right)$\textit{ has no zeros or poles on C. Then the total number
of zeros of }$F\left( s \right)$\textit{ in D, minus the total
number of poles of }$F\left( s \right)$\textit{ in D, counted with
multiplicities, is given by}
\[
\frac{1}{2\pi \,i}\int\limits_C {\frac{F'\left( s \right)}{F\left( s
\right)}} ds=\frac{1}{2\pi \,}\Delta_{C} \arg F\left( s \right).
\]

\vspace{2mm}

\textit{Here }$\Delta_{C} \arg F\left( s \right)$\textit{ denotes
the change of argument of }$F\left( s \right)$\textit{ along C.}

\vspace{3mm}

In addition we consider the following results obtained from the
Stirling's formula [16] and Jensen's formula ( [10], pp. 49-50):

\vspace{2mm}

\textbf{Proposition 6.1.} (Stirling's formula). \textit{We have}
\[
\log \Gamma \left( s \right)=\left( {s-\frac{1}{2}} \right)\,\log
s-s+\frac{1}{2}\log 2\pi +O\left( {\left| s \right|^{-1}}
\right)\approx \left( {s-\frac{1}{2}} \right)\,\log s-s+O\left( 1
\right)
\]
\textit{valid as} $\left| s \right|\to \infty $, \textit{in the
angle }$-\pi +\delta <\arg s<\pi -\delta $,\textit{ for any fixed
}$\delta
>0$.

\vspace{8mm}

\textbf{Proposition 6.2.} \textit{Let }$f$ \textit{be a function
which is analytic in a neighborhood of the disk }$\left| {z-a}
\right|<R$.

\textit{Suppose} $0<r<R$ \textit{and that f }has $n_{\, \, }$zeros
in the disk $\left| {z-a} \right|<r$. \textit{Let }$M=\max \left|
{f\left( {a+R\;e^{i\theta }} \right)\,} \right|$ \textit{and suppose
that }$\left| {f\left( 0 \right)} \right|\ne 0$. \textit{Then}.
\[
\left( {\frac{R}{r}} \right)^{n}\le \frac{M}{\left| {f\left( 0
\right)\,} \right|}.
\]

\vspace{3mm}

We begin considering the Theorem 6.1 for the function
$\textcyr{E1}\left( s \right)$ in the region $R$\textbf{, }whose $R$
is a rectangle in the complex plane with vertices at $2,\,2+iT,\,
-1+iT \,\mbox{and}\, -1$ (see Fig.1 and Appendix).

\vspace{3mm}

Let $D$ be the rectangular path passing through these vertices in
the anticlockwise direction.

\vspace{3mm}

It was noted earlier that $\textcyr{E1}\left( s \right)$ is analytic
everywhere, and has as its only zeros the imaginary zeros in the
critical strip.

\vspace{3mm}

Hence the number of zeros in the region $R$\textbf{,} which is given
by the equation

\vspace{2mm}

\[
N\left( T \right)=\frac{1}{2\pi \,i}\int\limits_D
{{\textcyr{E1}^{'}\left( s \right)} \mathord{\left/ {\vphantom
{{\textcyr{E1}^{'}\left( s \right)} \textcyr{E1} }} \right.
\kern-\nulldelimiterspace}\textcyr{E1} \left( s \right)\,\;ds=}
\frac{1}{2\pi \,}\Delta_{D} \arg \textcyr{E1} \left( s \right)
\]
and so
\[
2\pi \,N\left( T \right)=\Delta_{D} \arg \,\textcyr{E1}\left( s
\right).
\]

\vspace{2mm}

Our study of $N\left( T \right)$ will therefore focus on the change
of the argument of $\textcyr{E1}\left( s \right)$ as we move around
the rectangle $D$. As we move along the base of this rectangle,
there is no change in  $\arg \,\textcyr{E1}\left( s \right),$ since
$\textcyr{E1}\left( s \right)$ is real along this path and is never
equal to zero.

\vspace{3mm}

We wish to show that the change in $\arg \textcyr{E1} \left( s
\right)$\textbf{ }as $s$ goes from $\frac{1}{2}+iT$ to\textbf{
}$-1+iT$and then to $-1$ is equal to the change as $s$ moves from
$2$ to  $2+iT$ to $\frac{1}{2}+iT.$

\vspace{3mm}

To see this we observe that
\[
\textcyr{E1}\left( {\sigma +it} \right)=\textcyr{E1}\left( {1-\sigma
-it} \right)=\overline {\textcyr{E1}\left( {1-\sigma +it} \right)} .
\]

\vspace{1mm}

Hence the change in argument over the two paths will be the same.

\vspace{2mm}

If we define $L$ to be the path from $2$\textbf{ }to $2+iT$ then
$\frac{1}{2}+iT$\textbf{, }we have that

\vspace{2mm}

\[
2\pi \,N\left( T \right)=2\,\Delta_{L} \arg \textcyr{E1} \left( s
\right)
\]
or
\[
\pi \,N\left( T \right)=\Delta_{L} \arg \textcyr{E1}\left( s
\right).
\]

\vspace{3mm}

\medskip
\centerline{\epsfig{figure=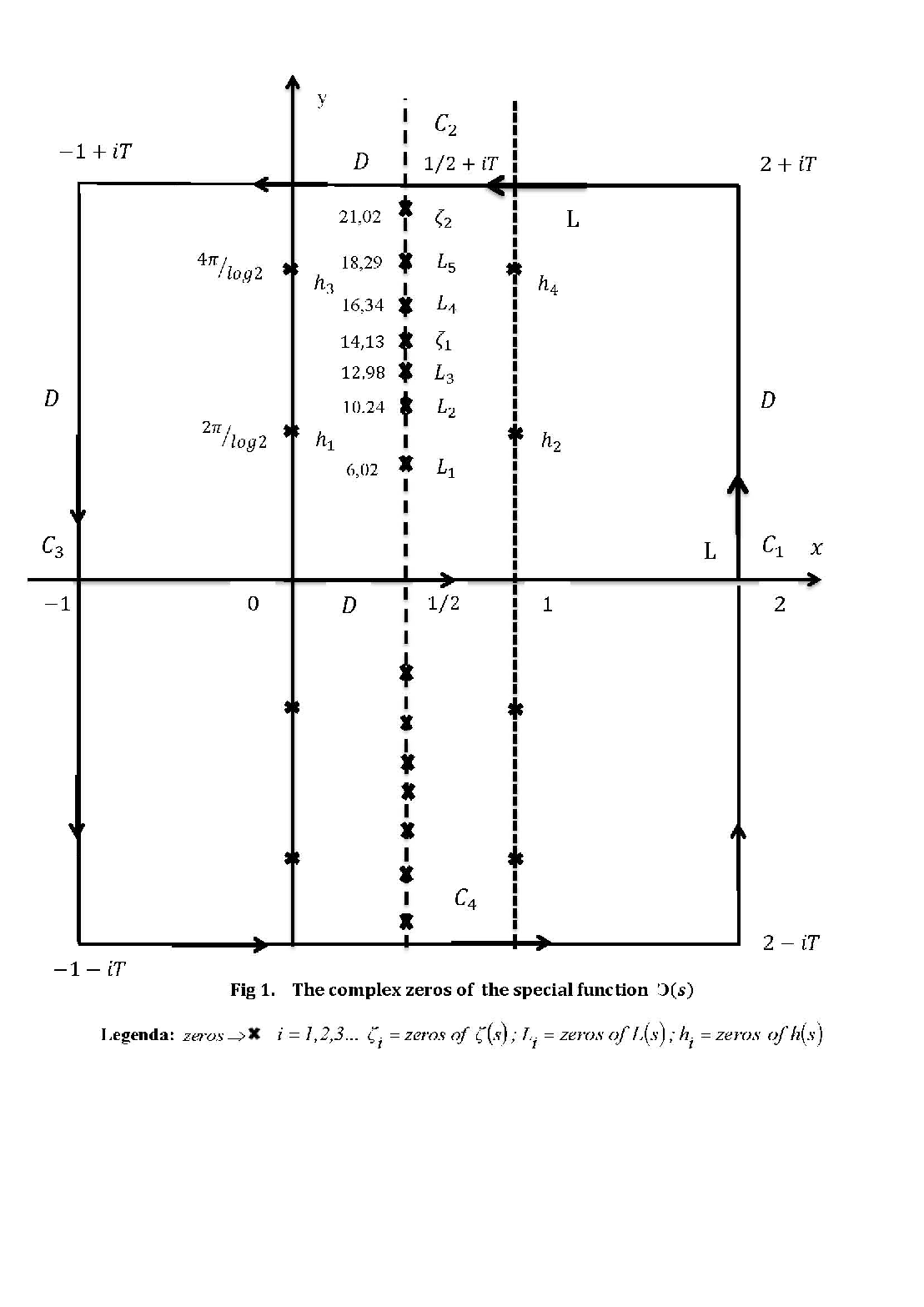,width=17.5cm,height=25cm}}
\medskip
\baselineskip=0.20in

We now recall the definition of $\textcyr{E1}\left( s \right)$ given
by
\[
\textcyr{E1}\left( s \right)=\frac{\left( {1-2^{s}} \right)\,\left(
{1-2^{1-s}} \right)\,\,\Gamma \left( s \right)\zeta \left( s
\right)L\left( s \right)}{\pi^{s}}
\]
and consider the argument of each section of the right-hand-size in
turn.

\vspace{2mm}

We have: $\Delta_{L} \arg \left[ {\left( {1-2^{s}} \right)\left(
{1-2^{1-s}} \right)} \right]=\,\Delta_{L} \arg \left( {1-2^{s}}
\right)+\Delta_{L} \arg \left( {1-2^{1-s}} \right)$
\[
\,=\,2\Delta_{L} \arg \left( {1-2^{1-s}} \right)=T\log 2+O\left( 1
\right)
\]
and
\[
\,\Delta_{L} \arg \pi^{-s}=\Delta_{L} \arg \exp \left( {-s\log \pi }
\right)=\,\Delta_{L} \left( {-t\log \pi } \right)=-T\log \pi .
\]

\vspace{1mm}

The proof of this first result is provided in Appendix.

\vspace{1mm}

To consider $\Gamma \left( s \right)$ we call on Stirling's formula
and also $\,\arg z=\Im \log z$, thus we have:
\[
\,\Delta_{L} \arg \Gamma \left( s \right)=\Im \log \Gamma \left(
{\frac{1}{2}+iT} \right)=\Im \left[ {iT\log \left( {\frac{1}{2}+iT}
\right)-\frac{1}{2}-iT+O\left( 1 \right)} \right]
\]

or since
\[
\,\log \left( {\frac{1}{2} + iT} \right) = \log \left| {\frac{1}{2}
+ iT} \right| + i\frac{\pi }{2} = \log \sqrt {\frac{1}{4} + T^2} +
i\frac{\pi }{2} \approx \log T + O\left( {\frac{1}{T}} \right) +
i\frac{\pi }{2}
\]
\[
\,\Delta_{L} \arg \Gamma \left( s \right)=T\log T-T+O\left( 1
\right).
\]

The above arguments can then be combined giving
\[
\pi \,N\left( T \right)=\Delta_{L} \arg \textcyr{E1} \left( s
\right)
\]
\[
=\Delta_{L} \arg \left[ {\left( {1-2^{s}} \right)\left( {1-2^{1-s}}
\right)} \right]+\Delta_{L} \arg \pi^{-s}+\Delta_{L} \arg \Gamma
\left( s \right)+\Delta_{L} \arg \zeta \left( s \right)+\Delta_{L}
\arg L\left( s \right)
\]
\[
=T\left( {\log 2-\log \pi +\log T-1} \right)+O(1)+\Delta_{L} \arg
\zeta \left( s \right)+\Delta_{L} \arg L\left( s \right)
\]
\[
=T\log \frac{2T}{\pi \,e}+\Delta_{L} \arg \zeta \left( s
\right)+\Delta _{L} \arg L\left( s \right)+O\left( 1 \right).
\]
Hence
\[
\,N\left( T \right)=\frac{T}{\pi }\log \frac{2T}{\pi \,e}+R\left( T
\right)+S\left( T \right)+O\left( 1 \right)
\]
where
\[
\pi \left[ {R\left( T \right)+S\left( T \right)} \right]=\Delta_{L}
\arg \zeta \left( s \right)+\Delta_{L} \arg L\left( s \right).
\]

From this point, in order to prove the approximation for $N\left( T
\right)$ initially claimed in (iii) it will be sufficient to show

\begin{equation}
\label{eq19} R\left( T \right)=S\left( T \right)=O\left( {\log T}
\right)\quad as \quad T\to \infty_{.}\end{equation}

\vspace{1mm}

First we need to know a bound for $\zeta \left( s \right)$ and
$L\left( s \right)$ on vertical strips.

\vspace{1mm}

Let $s=\sigma +it$ where $\sigma $ and $t$ are real.

\vspace{3mm}

\textbf{Proposition 6.3.} \textit{ Let }$0<\delta <1$. \textit{In
the region} $\sigma \ge \delta ,\;\,t>1$ \textit{we have}

\vspace{2mm}

$\quad \quad\quad \quad \quad\quad$(A) $\zeta \left( {\sigma +it}
\right)=O\left( {t^{1-\delta }} \right)$ ; (B) $L\left( {\sigma +it}
\right)=O\left( {t^{1-\delta }} \right)$.

\vspace{2mm}

\textbf{Proof.} Firstly we will deduce before the estimate (B). To
achieve this goal we use the following formula of Euler-Boole
summation\footnote{ NIST, Digital Library of Mathematical Functions,
(forthcoming) http://dlmf.nist.gov/24.17}, because it is used to
explain the properties of alternating series and it is better suited
than Euler-Maclaurin summation [4].

\vspace{2mm}

Let $0 \le h \le 1_{ }$ and $a,m$ and $n_{ }$ integers such $n > a$,
 $m > 0_{ }$ and $f^{\left( m \right)}\left( x \right)_{ }$ is
absolutely integrable over $\left[ {a,n} \right]_{.}$

\vspace{2mm}

Then we have:

\[\sum _{{j=a}}^{{n-1}}(-1)^{j}f(j+h)=\frac{1}{2}\sum _{{k=0}}^{{m-1}}\frac{\mathop{E_{{k}}\/}
\nolimits\!\left(h\right)}{k!}\left((-1)^{{n-1}}f^{{(k)}}(n)+(-1)^{a}f^{{(k)}}(a)\right)\]
\[
+\frac{1}{2(m-1)!}\int
_{a}^{n}f^{{(m)}}(x)\mathop{\widetilde{E}_{{m-1}}\/}\nolimits\!\left(h-x\right)dx.\]

$E_{n} \left( x \right)$ are Euler polynomials given by the
generating function:

\[
\frac{2e^{xt}}{e^{t}+1}=\sum\limits_{n=0}^\infty {E_{n} } \left( x
\right)\frac{t^{n}}{n\,\mbox{!}}.
\]

and the periodic Euler polynomials $\tilde {E}_n \left( x \right)$
are defined by setting  $\tilde {E}_n \left( x \right) = E_n \left(
x \right)$ for $0\le x\prec 1 $ and   $\tilde {E}_n \left( {x + 1}
\right) = - \tilde {E}_n \left( x \right)$ for all other x.

\vspace{2mm}

Let $N$ be a larger integer to be determined later.

\vspace{2mm}

If $f_{\, \, }$is any smooth function, for $M>N$, in the formula of
Euler-Boole summation above, with $a = N$, $m = 1$ and by taking the
limit as $h \to 0_{ }$ we obtain:

\[ \sum\limits_{n=N}^{M-1} {\left( {-1}
\right)^{n}} f\left( n \right)=\frac{1}{2}\left[ {\left( {-1}
\right)^{N}f\left( N \right)+\left( {-1} \right)^{M-1}f\left( M
\right)} \right]+\frac{1}{2}\int\limits_N^M \tilde {E}_0 \left(-x
\right)f'\left( x \right)dx
\]

\vspace{1mm}

where $\tilde E_0 \left( x \right) = sgn\left[ {\sin \left( {\pi
\,x} \right)} \right]$, that is a pieciewise constant periodic
function.

\vspace{1mm}

Take $f\left( x \right)=\left( {2x+1} \right)^{-s}$, where initially
$\Re \left( s \right)>1$, and let $M\to \infty $.

\vspace{1mm}

We obtain:
\[
L\left( s \right)-\sum\limits_{n\prec N} {\,\frac{\left( {-1}
\right)^{n}}{\left( {2n+1} \right)^{s}}=} \sum\limits_{n=N}^\infty
{\frac{\left( {-1} \right)^{n}}{\,\left( {2n+1}
\right)^{s}}=\frac{1}{2}\left[ {\left( {-1} \right)^{N}\left( {2N+1}
\right)^{-s}} \right]} -s\int\limits_N^\infty  \tilde {E}_0 \left(-
x \right)\left( {2x+1} \right)^{-s-1}dx.
\]

The integral $s\int\limits_N^\infty  \tilde {E}_0 \left(- x
\right)\left( {2x+1} \right)^{-s-1}dx$ is absolutely convergent if
$\sigma =\Re \left( s \right)>0$, and since $\left|{\tilde E_0
\left( { - x} \right)} \right| = 1$, we note that

\[
\left| {s\int\limits_N^\infty  \left( {2x+1} \right)^{-s-1}dx}
\right|<\left| s \right|\int\limits_N^\infty \left( {2x+1}
\right)^{-\sigma -1}dx=\frac{\left| s \right|}{\sigma }\left( {2N+1}
\right)^{-\sigma }\le \left( {1+\frac{t}{\sigma }} \right)\,\left(
{2N+1} \right)^{-\sigma }
\]
where we have used the triangle inequality $\left| s \right|\le
\sigma +t$.

Also
\[
\left| {\sum\limits_{n\prec N}{\frac{\left( {-1} \right)^{n}}{\left(
{2n+1} \right)^{s}}} } \right|\le \sum\limits_{n\prec N}
{\frac{\left( {-1} \right)^{n}}{\left( {2n+1} \right)^{\sigma }}<}
\int\limits_0^N \left( {2x+1} \right)^{-\sigma }dx=\frac{\left(
{2N+1} \right)^{1-\sigma }}{1-\sigma }-\frac{1}{1-\sigma }.
\]

Thus
\begin{equation}
\label{eq20}\ \left| {L\left( s \right)} \right|=\left|
{\sum\limits_{n\prec N} {\frac{\left( {-1} \right)^{n}}{\left(
{2n+1} \right)^{s}}+\frac{1}{2}\left[ {\left( {-1} \right)^{N}\left(
{2N+1} \right)^{-s}} \right]-s\int\limits_N^\infty  \left( {2x+1}
\right)^{-s-1}dx} } \right|\end{equation}

\vspace{1mm}

\[\le \frac{\left( {2N+1}
\right)^{1-\sigma }}{1-\sigma }-\frac{1}{1-\sigma
}+\frac{1}{2}\left[ {\left( {-1} \right)^{N}\left( {2N+1}
\right)^{-\sigma }} \right]+\left( {1+\frac{t}{\sigma }}
\right)\,\left( {2N+1} \right)^{-\sigma }\]

\vspace{1mm}

\[ <\frac{\left( {2N+1} \right)^{1-\sigma }}{1-\sigma
}+\left( {\frac{3}{2}+\frac{t}{\sigma }} \right)\,\left( {2N+1}
\right)^{-\sigma }.\]

\vspace{1mm}

Assuming that $t>1$, we may estimate this by taking $N$ to be
greatest integer less than $\left( {\frac{t-1}{2}} \right).$

To see that this is the optimal choice of $t$, consider the two
potentially largest terms in (\ref{eq20}):

\[\frac{\left( {2N+1} \right)^{1-\sigma }}{1-\sigma }\quad  and \quad
\left( {\frac{t}{\sigma }} \right)\,\left( {2N+1} \right)^{-\sigma
}.\]

If we take $N$ to be approximately $\frac{t^{\alpha }-1}{2}$ for
some $\alpha $, these are $\frac{t^{\alpha \,\left( {1-\sigma }
\right)}}{\left( {1-\sigma } \right)}$ and $\left( \sigma \right)^{
- 1} \cdot t^{1 - \alpha \sigma }.$

\vspace{1mm}

As $\alpha $ varies, one increases, the other decreases.

\vspace{1mm}

Thus, we want to equate the exponents, so $ \alpha \left( {1-\sigma
} \right)=1-\alpha \sigma $, or $\alpha =1$.

\vspace{1mm}

Taking $N\approx \frac{t-1}{2}$ , we see that $L\left( s \right)$ is
of the order $O\left( {t^{1-\sigma }} \right)$.

\vspace{2mm}

If $\sigma \ge \delta $ and $t>1$, we see that $L\left( {\sigma +it}
\right)=O\left( {t^{1-\delta }} \right)_{\, \, }$and thus (B) is
proved.$_{\, }$

\vspace{2mm}

To achieve the estimate (A) it is sufficient to use the same
procedure, but in this case we recall the formula of
Euler-Maclaurin, that is
\[
\sum\limits_{n=N}^M  f\left( n \right)=\int\limits_N^M {f\left( x
\right)} dx+\frac{1}{2}f\left( N \right)+\frac{1}{2}f\left( M
\right)\;+\;\int\limits_N^M {B_{1} } \left( {x-\left[ x \right]}
\right)f'\left( x \right)dx
\]
where $B_{1} \left( x \right)=x-\frac{1}{2}$ is the first Bernouilli
polynomial, $\left[ x \right]$ is the greatest integer and take
$f\left( x \right)=x^{-s}$, where initially $\Re \left( s
\right)>1$, and let $M\to \infty $.

\vspace{2mm}

In this case, at the end, we obtain

\[
\left| {\zeta \left( s \right)} \right|\le \frac{\left( N
\right)^{1-\sigma }}{1-\sigma }+\frac{\left( N \right)^{1-\sigma
}}{t}+\left( {\frac{1}{2}+\frac{t}{2\sigma }} \right)\,\left( N
\right)^{-\sigma }.
\]

Taking $N\approx t_{\, }$, we see that, if $\sigma \ge \delta $ and
$t>1$, $\zeta \left( {\sigma +it} \right)=O\left( {t^{1-\delta }}
\right)_{.}$

\vspace{1mm}

Consequently (A) is proved.

\vspace{2mm}

Finally we will prove (\ref{eq19}), that is the integrals
\[
\Im \;\left( {\int\limits_2^{\frac{1}{2}+iT} {\frac{\zeta '\left( s
\right)}{\zeta \left( s \right)}ds+\int\limits_2^{\frac{1}{2}+iT}
{\frac{L'\left( s \right)}{L\left( s \right)}ds} } } \right)=O\left(
{\log T} \right).
\]

Firstly we note that $\zeta \left( s \right)$ and $L\left( s
\right)$ are holomorphic and non-vanishing in the half plane $\Re
\left( s \right)\succ 1$.

\vspace{1mm}

If $T$ is real, we have
\[\int\limits_2^{2+iT} {\frac{\zeta '\left( s \right)}{\zeta \left(
s \right)}ds} =\log \zeta \left( {2+iT} \right)-\log \zeta \left( 2
\right)\] and
\[
\int\limits_2^{2+iT} {\frac{L'\left( s \right)}{L\left( s
\right)}ds} =\log L\left( {2+iT} \right)-\log L\left( 2 \right).
\]
Here

\[\left| {\zeta \left( {2+iT} \right)-1} \right|=\left|
{\sum\limits_{n=2}^\infty {n^{-2-it}} } \right|\le
\sum\limits_{n=2}^\infty {\left| {n^{-s}} \right|} =\zeta \left( 2
\right)-1=0.644934\]
 and
\vspace{1mm}
\[
\left| {L\left( {2+iT} \right)-1} \right|=\left|
{\sum\limits_{n=1}^\infty {\left( {-1} \right)^{n}\left( {2n+1}
\right)^{-2-it}} } \right|<\left| {\,\sum\limits_{n=2}^\infty \left(
{2n-1} \right)^{-2-it}} \right|
\]
\[
<\left| {\left( {1-2^{-s}} \right)}
\right|\,\sum\limits_{n=2}^\infty {\left| {n^{-s}} \right|}
=\frac{3}{4}\left( {\zeta \left( 2 \right)-1} \right)=0.4837
\]

Since these are less than $1$, $\zeta \left( {2+iT} \right)$ and
$L\left( {2+iT} \right)$ are constrained to a circle which excludes
the origin, and
\begin{equation}
\label{eq21}\left| {\zeta \left( {2+iT} \right)} \right|>1-0.644934
\quad and \quad \left| {L\left( {2+iT} \right)}
\right|>1-0.4837\end{equation}

Finally, we have that

\begin{equation}
\label{eq22} \int\limits_2^{2+iT} {\frac{\zeta '\left( s
\right)}{\zeta \left( s \right)}ds} =O\left( 1 \right)_{\, \, }\quad
and \quad \int\limits_2^{2+iT} {\frac{L'\left( s \right)}{L\left( s
\right)}ds} =O\left( 1 \right)\end{equation}

\vspace{1mm}

To complete the proof of (\ref{eq19}) we show that

\[\Im \;\left(
{\int\limits_{2+iT}^{\frac{1}{2}+iT} {\frac{\zeta '\left( s
\right)}{\zeta \left( s
\right)}ds+\int\limits_{2+iT}^{\frac{1}{2}+iT} {\frac{L'\left( s
\right)}{L\left( s \right)}ds} } } \right)=O\left( {\log T}
\right).\]

\vspace{1mm}

We assume that the path from $2+iT$ to $\frac{1}{2}+iT$ does not
pass through any zero of $\zeta \left( s \right)$ and any zero of
$L\left( s \right)$, by moving the path up slightly if necessary.

\vspace{1mm}

By the Argument Principle the two integrals represent respectively
the change in the argument of $\zeta \left( s \right)$ and the
change in the argument of $L\left( s \right)$ as $s$ moves from
$2+iT$ to $\frac{1}{2}+iT$.

\vspace{1mm}

These are approximately $\pi \left( {c_{1} +c_{2} } \right)$, where
$c_{1} $ is the number of sign changes in $\Re \zeta \left( {s+it}
\right)$ and $c_{2} $ is the number of sign changes in $\Re L\left(
{s+it} \right)$, as $s_{\, }$ moves from $2$ to $\frac{1}{2}$, since
the sign must change every time the argument changes by $\pi $.

\vspace{1mm}

We note that if $s$ is real:$\qquad\Re \zeta \left( {s+it}
\right)=\frac{1}{2}\left[ {\zeta \left( {s+it} \right)+\zeta \left(
{s-it} \right)} \right]$ and $ \Re L\left( {s+it}
\right)=\frac{1}{2}\left[ {L\left( {s+it} \right)+L\left( {s-it}
\right)} \right]  $.

\vspace{1mm}

Therefore, it is sufficient to show that the number of zeros of
$\frac{1}{2}\left[ {\zeta \left( {s+it} \right)+\zeta \left( {s-it}
\right)} \right]$ and the number of zeros of $\frac{1}{2}\left[
{L\left( {s+it} \right)+L\left( {s-it} \right)} \right]$ on the
segment $\left[ {\frac{1}{2},2} \right]$ of real axis are $O\left(
{\log T} \right)$.

\vspace{1mm}

In fact, we will use Proposition 6.2 to estimate the number of zeros
of $f\left( s \right)=\frac{1}{2}\left[ {\zeta \left( {s+it}
\right)+\zeta \left( {s-it} \right)} \right]$ and the number of
zeros of $g\left( s \right)=\frac{1}{2}\left[ {L\left( {s+it}
\right)+L\left( {s-it} \right)} \right]$ inside the circle $\left|
{s-2} \right|\prec \frac{3}{2}$.

\vspace{1mm}

We take $a=2,\;R=\frac{7}{4}\;\mbox{and}\;r=\frac{3}{2}$ in the
Proposition 6.2.

\vspace{1mm}

First, we note that $\left|f\left( 2 \right)\right|$ and
$\left|g\left( 2 \right)\right|$ are bounded by (\ref{eq21}).

On the other hand,
\[
\mathop {\max }\limits_{\left| {s - 2} \right| = 7 \mathord{\left/
{\vphantom {7 4}} \right. \kern-\nulldelimiterspace} 4} \left|
{f\left( s \right)} \right| = O\left( {T^{3 \mathord{\left/
{\vphantom {3 4}} \right. \kern-\nulldelimiterspace} 4}}
\right)\quad \;and\quad \;\mathop {\max }\limits_{\left| {s - 2}
\right| = 7 \mathord{\left/ {\vphantom {7 4}} \right.
\kern-\nulldelimiterspace} 4} \left| {g\left( s \right)} \right| =
O\left( {T^{3 \mathord{\left/ {\vphantom {3 4}} \right.
\kern-\nulldelimiterspace} 4}} \right)
\]

by Proposition 6.3.

Therefore if $_{\, }n$ is the number of zeros of $f\left( s \right)$
inside $\left| {s-2} \right|<\frac{3}{2}$ and if $m$ is the number
of zeros of $g\left( s \right)$ inside $\left| {s-2}
\right|<\frac{3}{2}$ , we have
\[\left( {\frac{7 \mathord{\left/ {\vphantom {7 4}} \right.
\kern-\nulldelimiterspace} 4}{3 \mathord{\left/ {\vphantom {3 2}}
\right. \kern-\nulldelimiterspace} 2}} \right)^{n}=O\left( {T^{3
\mathord{\left/ {\vphantom {3 4}} \right. \kern-\nulldelimiterspace}
4}} \right)\quad \;and\quad\left( {\frac{7 \mathord{\left/
{\vphantom {7 4}} \right. \kern-\nulldelimiterspace} 4}{3
\mathord{\left/ {\vphantom {3 2}} \right. \kern-\nulldelimiterspace}
2}} \right)^{m}=O\left( {T^{3 \mathord{\left/ {\vphantom {3 4}}
\right. \kern-\nulldelimiterspace} 4}} \right)\,,\] or taking
logarithms in the first case we have that $n\log \left( {7
\mathord{\left/ {\vphantom {7 6}} \right. \kern-\nulldelimiterspace}
6} \right)$ is bounded by $\frac{3}{4}\log \left( T \right)$ plus a
constant and the latter case we have that $m\log \left( {7
\mathord{\left/ {\vphantom {7 6}} \right. \kern-\nulldelimiterspace}
6} \right)$ is bounded by $\frac{3}{4}\log \left( T \right)$ plus a
constant.

This completes the proof of (iii).

\vspace{1mm}
\section{Conclusion}

The symmetric form of the functional equation for $\zeta \left( s
\right)$ represents one of the most remarkable achievements by B.
Riemann.

This fundamental result was discovered and proved in his paper [15]
in two different ways: the first was described in Section 2, the
latter is similar to the one that we have illustrated in Section 3:
it is conceptually more difficult because required taking the Mellin
transform to boot and use an integral involving the theta function.

All modern proofs of the functional equation involve mathematical
tools that were unavailable to L. Euler and it is remarkable that he
was nevertheless able to predict the asymmetric form of the
functional equation for the Zeta function.

In his paper [8] Euler used the differentiation of divergent series
and a version of his of the Euler-Maclaurin summation formula.

Here we presented a proof of symmetric form of the functional
equation for the Zeta function that Euler himself could have proved
with a little step at end of his paper.

The result of this simple proof, based upon the three reflection
formulae of $\eta \left( s \right)$, $L\left( s \right)$ and $\Gamma
\left( s \right)$ with the duplication formula of \textit{sine}, is
a most general form of the functional equation for Riemann Zeta
function.

It is easy to see that if $f\left( s \right)$ and $g\left( s
\right)$ are two Dirichlet series, each satisfying a functional
equation, then the product $f\left( s \right)\cdot g\left( s
\right)$ defines a third Dirichlet series also satisfying a given
functional equation, but, in our specific case, with the product of
two functional equation in the asymmetric form we have obtained a
functional equation in the unexpected symmetric form.

In the first part of this paper we obtained also an amazing
connection between the Jacobi's imaginary transformation and an
infinite series identity of Ramanujan.

In the second part using techniques similar to those of Riemann, it
is shown how to locate and count the imaginary zeros of the entire
function $\textcyr{E1} \left( s \right)$, which is an extension of
the special function $A\left( s \right)$, that we have previously
introduced [14].

Here we apply also the Euler-Boole summation formula and we obtain
an estimate of the distribution of the zeros of the function
$\textcyr{E1} \left( s \right)$ to follow a method, which Ingham
([10] pp.68-71) attributes to Backlund [1].

Basically we use the fact that we have a bound on the growth of
$\zeta \left( s \right)$ and the growth of $L\left( s \right)$ in
the critical strip.

More precisely with the Fundamental Theorem we also established that
the number of the zeros of the function $\textcyr{E1} \left( s
\right)$ in the critical strip is :

\begin{equation}
\label{eq23}N_{\textcyr{E1}} \left( T \right)=\frac{T}{\pi }\log
\frac{2T}{\pi \, e}+O\left( {\log T} \right)\end{equation}

Now, from Appendix, we have that the number of zeros of the factor
$h\left( s \right)=\left( {1-2^{s}} \right)\,\left( {1-2^{1-s}}
\right)$ is :
\begin{equation}
\label{eq24} N_{h} \left( T \right)=\frac{T}{\pi }\log 2+O\left( 1
\right).\end{equation}

Subtracting (\ref{eq24}) from (\ref{eq23}) we have the number of
zeros of the special function $A\left( s \right)$, that is:
\begin{equation}
\label{eq25} N_{A} \left( T \right)=\frac{T}{\pi }\log \frac{T}{\pi
\, e}+O\left( {\log T} \right).\end{equation}

 and from ( [18],p. 214, 9.4.3 ) we have that the distribution function for the zeros of the
Riemann Zeta function is :
\begin{equation}
\label{eq26} N_{\zeta } \left( T \right)=\frac{T}{2\pi }\log
\frac{T}{2\pi \,}-\frac{T}{2\pi }+O\left( {\log T}
\right)=\frac{T}{2\pi }\log \frac{T}{2\pi \,  e}+O\left( {\log T}
\right).\end{equation}
Now subtracting (\ref{eq26}) from(\ref{eq25})
we have the number of zeros of the Dirichlet $L$ function:
\begin{equation}
\label{eq27} N_{L} \left( T \right)=\frac{T}{2\pi }\log
\frac{2T}{\pi \, e}+O\left( {\log T} \right).\end{equation} The
previous results describe, in detail, the structure of the complex
roots of the entire function $\textcyr{E1}\left( s \right)$.

Table 1 (see the following page) shows the frequency distribution
for the actual zeros in successive intervals of  $t$ .

The author used  M. Rubinstein's $L$-function calculator\footnote{
http://oto.math.uwaterloo.ca/$\sim
$mrubinst/L{\_}function{\_}public/L.html\par } to compute, with
approximation, the complex zeros in the critical line $\sigma = 1
\mathord{\left/ {\vphantom {1 2}} \right. \kern-\nulldelimiterspace}
2\;$ and in the interval $\;0 \le t \le 100\;$ (see Fig. 2).

\vspace{3mm}

\medskip
\centerline{\epsfig{figure=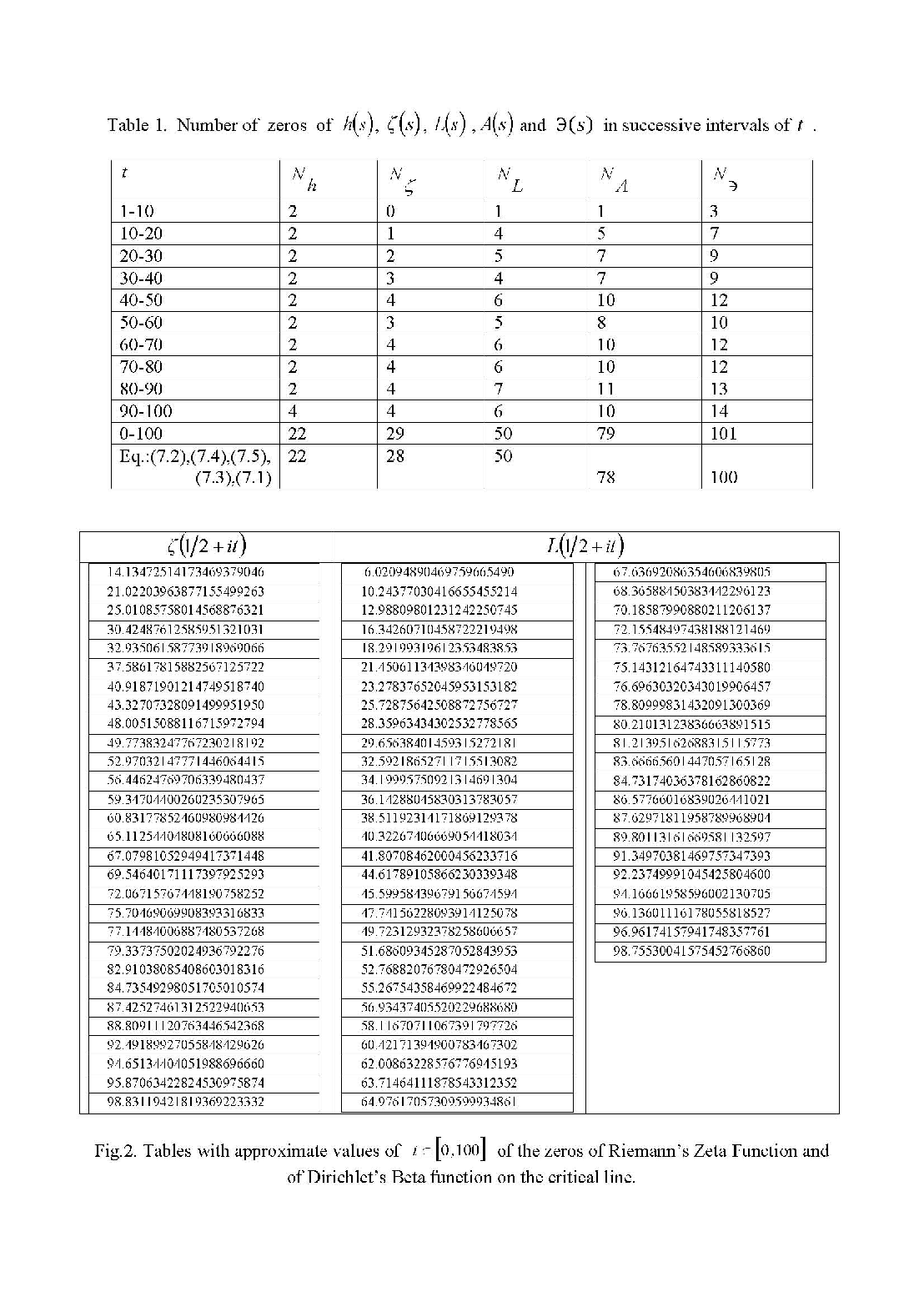,width=16cm,height=23cm}}
\medskip
\baselineskip=0.20in

\vspace{1mm}

This makes the strong difference in the distributions of the gaps,
all very interesting.

\vspace{1mm}

In this case, having also $N_{A} \ge N_{A} \left( {T=100} \right)$,
it follows that there are exactly $N$ zeros in this portion of the
critical strip, all lying on the critical line.

\vspace{1mm}

To be complete, we give also, for large $T$, the following result:
\[N_{L} \left( T \right)=N_{\zeta } \left( T \right)+N_{h} \left( T
\right) \quad and \quad N_{\textcyr{E1}} \left( T \right)=2\,N_{L}
\left( T \right)\;.\]

In addition we observe that the complex roots of the factor $h\left(
s \right)$ lie on the vertical lines $\Re \left( s \right)=0$ and
$\Re \left( s \right)=1$ and they are separated by $\frac{2\pi
\;i}{\log 2}$.

\vspace{1mm}

While if we assume the Generalized Riemann Hypothesis
(GRH)\footnote{ GRH: Riemann Hypothesis is true and in addition the
nontrivial zeros of all Dirichlet $L$-functions lie on the critical
line $\Re \left( s \right)=1 \mathord{\left/ {\vphantom {1 2}}
\right. \kern-\nulldelimiterspace} 2\;.$}, this implies that all
complex zeros of the special function $A\left( s \right)$ lie on the
vertical line $\Re \left( s \right)=\frac{1}{2}$ and thus, at a
height $T$ the average spacing between zeros is asymptotic to
$\frac{\pi }{\log T}$.

\section{Appendix}

\vspace{5mm}

We study the solution in $s$ of the following Dirichlet polynomial:

\begin{equation}
\label{eq28} f\left( s \right)=1-2^{1-s}=1-2\left( {\frac{1}{2}}
\right)^{s}=0.\end{equation}

\vspace{3mm}

This is the simplest example of a Dirichlet polynomial equation.

\vspace{3mm}

In this case, the complex roots are

\[s=1\pm \frac{2\pi \,i\,k}{\log 2} \quad  with \quad k\in Z_{.}\]

Hence the complex roots lie on the vertical line $\Re \left( s
\right)=1$ and are separated by $\frac{2\pi \;i}{\log 2}$.

\vspace{3mm}

In order to establish the density estimate of (\ref{eq28}), we will
estimate the winding number of the function $f(s)=1-2\left(
{\frac{1}{2}} \right)^{s}$ when $s$ runs around the contour $C_{1}
+C_{2} +C_{3} +C_{4} $, where $C_{1} $  and   $C_{3} $ are the
vertical line segments $2-iT\to 2+iT$ and $-1+iT\to -1-iT$ and
$C_{2} $ and $C_{4} $ are the horizontal line segments $2+iT\to
-1+iT$ and $-1-iT\to 2-iT$, with $T>0$ (see Fig.1).

\vspace{1mm}

For $\Re \left( s \right)=2$ we have$\left| {1-f\left( s \right)}
\right|=\left| {2\left( {\frac{1}{2}} \right)^{s}}
\right|=\frac{1}{2}<1$, so the winding number along $\;C_{1}\; $is
at most $\frac{1}{2}$.

\vspace{1mm}

Likewise, for $\Re \left( s \right)=-1$, we have $1<\left| {f\left(
s \right)-1} \right|\,=\left| {2\left( {\frac{1}{2}}
\right)^{-1+iT}} \right|\le 2\left( {\frac{1}{2}} \right)^{-1}=4$ so
the winding number along $C_{3} $ is that of term $2\left(
{\frac{1}{2}} \right)^{s}$, up to at most $\frac{1}{2}$.

\vspace{1mm}

Hence, the winding number along the contour $C_{1}+C_{3} $ is equal
to $\left( {\frac{T}{\pi }} \right)\log 2$, up to at most $1$.

\vspace{1mm}

We will now show that the winding number along $C_{2} +C_{4} $ is
bounded, using a classical argument ([10], p. 69).

\vspace{1mm}

Let $n$ the number of distinct points on $C_{2} $ at which $\Re
\,f\left( s \right)=0$.

\vspace{1mm}

For real value of $z$,

\[
\Re \,f\left( {z+iT} \right)=\frac{1}{2}\left[ {f\left( {z+iT}
\right)+f\left( {z-iT} \right)} \right].
\]

\vspace{3mm}

Hence, putting $g\left( z \right)=2\,\Re \,f\left( {z+iT} \right)$
we see that $n$ is bounded by the number of zeros of $g$ in a disk
containing the interval $\left( {0,1} \right)$.

\vspace{3mm}

We take the disk centred at $2$, with radius $3$.

\vspace{3mm}

We have
\[
\left| {g\left( 2 \right)} \right|\ge 2-2\cdot 2\left( {\frac{1}{2}}
\right)^{2}=1>0.
\]

Furthermore, let $G$ the maximum of $g$ on disk with the same centre
and radius $\mbox{\it e}\cdot \left( 3 \right)$, so
\[
G\le 2+2\cdot \left( {\frac{1}{2}} \right)^{2-\mbox{\it e}\cdot
3}\;.
\]

By Proposition 6.2, it follows that $n\le \log \,\left| {G
\mathord{\left/ {\vphantom {G {g\left( 2 \right)}}} \right.
\kern-\nulldelimiterspace} {g\left( 2 \right)}} \right|.$

\vspace{1mm}

This gives a uniform bound on the winding number over $C_{2} $. $ $
The winding number over $C_{4} $ is estimated in the same manner.

\vspace{2mm}

We conclude from the above discussion that the winding number of
$f\left( s \right)=1-2^{1-s}$ over the closed contour $C_{1} +C_{2}
+C_{3} +C_{4} $ equals $\left( {\frac{T}{\pi }} \right)\log 2$, up
to a constant (dependent on $f)$, from which follows the asymptotic
density estimate:
\[
D_{f} =\left( {\frac{T}{\pi }} \right)\log 2+O\left( 1 \right).
\]
If we count the zeros in the upper half of a vertical strip $\left\{
{s\mbox{:\, 0}\le \Im \left( s \right)\le T} \right\}$ we have:
\[
N_{f} =\left( {\frac{T}{2\,\pi }} \right)\log 2+O\left( 1 \right).
\]

A lot of details in relation to what we have just shown were
published in ([12], chap. 3, pp. 63-77).

\section{Additional Remark}

\vspace{2mm}

The author is aware that some of the results presented in [14] and
in this paper are not new.

In particular, the main subject of this paper, the function
$\textcyr{E1}\left( s \right)$ is, apart from a factor $\left(
{1-2^{s}} \right)\,\left( {1-2^{1-s}} \right)\,\;{\Gamma \left( s
\right)} \mathord{\left/ {\vphantom {{\Gamma \left( s \right)}
{\pi^{s}}}} \right. \kern-\nulldelimiterspace} {\pi^{s}},$ equal to
the product of the Riemann Zeta function and a certain $L$-function.

That product is equal to the Dedekind Zeta function associated to
the algebraic number field obtained from the field of rational
number by adjoining a square root of -1.

Let $r_{2} \left( n \right)$ denote the number of ways to write $n$
as sum of two squares, then the generating series for $r_{2} \left(
n \right)$:
\[
\zeta _{Q\left( {\sqrt { - 1} } \right)} \left( s \right) =
\frac{1}{4}\sum\limits_{n = 1}^\infty {r_2 \left( n \right)\,}
\left( n \right)^{ - s}
\]
is precisely the Dedekind Zeta function of the number field $Q\left(
{\sqrt { - 1} } \right)$, because it counts the number of ideals of
norm n.

It factors as the product of two Dirichlet series:
\[
\zeta _{Q\left( {\sqrt { - 1} } \right)} \left( s \right) = \zeta
\left( s \right)\,L\left( {s,\chi _4 } \right).
\]
The factorization is a result from class field theory, which
reflects the fact that an odd prime can be expressed as the sum of
two squares if and only if it is congruent to 1 modulo 4.

Dedekind Zeta functions were invented in the 19th century, and in
the course of time many of their properties have been established.
Some of the present results are therefore special cases of
well-known properties of the Dedekind Zeta.

Nevertheless, the goal of this manuscript is to highlight some
demonstration, direct and by increments, for treating certain
functional equations and special functions involved, as inspired by
methods similar to the ones used by Euler in his paper [8] and in
many other occasions (see [20], [19], and [21], chap. 3).

In order not to leave unsatisfied the reader's curiosity, we recall
that the choice of the letter \textcyr{E1} for the special function
\textcyr{E1}(s)\ is in honour of \textcyr{E1i0ler} (Euler).

\vspace{3mm}

\noindent Address: Via delle Azzorre 352/D2, Roma - Italy
\vspace{1mm}

\noindent E-mail: andrea.ossicini@yahoo.it

\end{document}